\documentclass[11pt,a4paper]{amsart}
\usepackage{amsmath,amssymb,mathtools}
\usepackage[margin=27mm]{geometry}
\usepackage{microtype}
\usepackage[hidelinks]{hyperref}
\hypersetup{
  pdftitle={Branched Covers and Fiberwise Degeneration along Weak K\"ahler Geodesics},
  pdfauthor={Guo Wenjiong and Sheng Zekun}
}
\numberwithin{equation}{section}
\newtheorem{theorem}{Theorem}[section]
\newtheorem{proposition}[theorem]{Proposition}
\newtheorem{lemma}[theorem]{Lemma}
\newtheorem{corollary}[theorem]{Corollary}
\theoremstyle{remark}
\newtheorem{remark}[theorem]{Remark}
\newcommand{\ddc}{dd^c}
\newcommand{\ii}{\sqrt{-1}}
\newcommand{\PP}{\mathbb P}
\newcommand{\FS}{\mathrm{FS}}
\newcommand{\pr}{\operatorname{pr}}
\newcommand{\supp}{\operatorname{supp}}

\newcommand{\HH}{\mathcal H}
\newcommand{\RR}{\mathcal R}
\newcommand{\norm}[1]{\left\lVert#1\right\rVert}

\newcommand{\wh}[1]{\widehat{#1}}
\allowdisplaybreaks[1]
\title[Fiberwise degeneration along weak K\"ahler geodesics]{Branched Covers and Fiberwise Degeneration along Weak K\"ahler Geodesics}
\author[Guo Wenjiong]{Guo Wenjiong\textsuperscript{1}}
\author[Sheng Zekun]{Sheng Zekun\textsuperscript{2}}
\thanks{\textsuperscript{1}Department of Mathematics, Fudan University. E-mail: \href{mailto:26110180012@m.fudan.edu.cn}{\texttt{26110180012@m.fudan.edu.cn}}.}
\thanks{\textsuperscript{2}Department of Mathematics, Nanjing University. E-mail: \href{mailto:221840029@smail.nju.edu.cn}{\texttt{221840029@smail.nju.edu.cn}}.}
\date{November 5, 2026}

\begin{document}
\begin{abstract}
On the projective line, we construct a pair of smooth K\"ahler potentials whose weak geodesic has isolated quadratic zeros of the fiber form at two fixed points throughout an interior time interval. The geodesic is smooth near these zeros, and its space-time form has constant rank one there. For the same fixed endpoints, the minimum fiber eigenvalue of the nondegenerate Monge--Amp\`ere approximations tends to zero, answering AIMPL Problem~2.5 in the negative. The construction works on every compact connected Riemann surface with a prescribed K\"ahler form and realizes any positive even order of vanishing. Taking products yields higher-dimensional examples in which the fiber rank drops by one along a smooth hypersurface while the space-time form retains maximal rank nearby. The proof starts from a semipositive model built from a branched cover. We make its endpoints strictly positive by local corrections and use comparison along holomorphic curves to prove that the weak geodesic still agrees with the model on an explicit interior region. There a ramification point of local degree $m$ gives a fiber zero of order exactly $2m-2$.
\end{abstract}
\maketitle

\section{Introduction and main results}

Along a weak K\"ahler geodesic, the rank of the space-time form and the rank of its restriction to a spatial fiber record different information. We construct segments with smooth strictly K\"ahler endpoints for which the space-time form has maximal rank near a ramification point, yet its restriction to the fiber vanishes at that point throughout an interior time interval. We first give explicit endpoints on $\PP^1$ and compute the quadratic zeros of the fiber form. We then carry out the construction on an arbitrary compact Riemann surface with prescribed K\"ahler form.

Let $(X,\omega)$ be a compact K\"ahler manifold of complex dimension $n$, and write
\[
 \HH_\omega=\{u\in C^\infty(X,\mathbb R):\omega+\ddc u>0\},
 \qquad \ddc=\ii\partial\bar\partial.
\]
Set $A=\{\tau\in\mathbb C:1<|\tau|<e\}$ and $t=\log|\tau|$. The weak geodesic joining $u_0,u_1\in\HH_\omega$ is represented by the $S^1$-invariant function $U(x,\tau)=u_t(x)$ satisfying
\begin{equation}\label{eq:hcma}
 \Omega_U:=\pr_X^*\omega+\ddc U\ge0,
 \qquad \Omega_U^{n+1}=0,
 \qquad U|_{|\tau|=e^j}=u_j\quad(j=0,1).
\end{equation}
The equation is understood in the Bedford--Taylor sense \cite[Proposition~2.3]{bedford-taylor}. For smooth K\"ahler endpoints, existence follows from Chen \cite[Theorem~3]{chen}. By \cite[Lemma~2.1 and equation~(5)]{darvas-survey}, $U$ is the unique bounded solution of \eqref{eq:hcma} and is the upper semicontinuous regularization of the supremum of all $S^1$-invariant $\pr_X^*\omega$-plurisubharmonic functions with these boundary upper bounds. The solution has real $C^{1,1}$ regularity \cite[Theorem~1.1]{ctw}. We call $\Omega_U$ the \emph{total form} and $\omega+\ddc_Xu_t$ the \emph{fiber form}.

Darvas's Question~2.3 \cite{darvas-survey} asks whether every such segment admits $\varepsilon>0$, depending on the fixed background and endpoints, such that
\begin{equation}\label{eq:darvas-question}
 \omega+\ddc_Xu_t\ge\varepsilon\omega\qquad(0\le t\le1)
\end{equation}
in the sense of currents. Xiao and Yao give a negative answer on every compact K\"ahler manifold \cite[Theorem~A]{xiao-yao}. Their weak geodesic is independent of time on an open region where the fiber rank is $n-1$, so the total form there also has rank $n-1$. For the explicit endpoints below, the fiber form has isolated quadratic zeros and the total form retains the maximal rank allowed by \eqref{eq:hcma} throughout a neighborhood of those zeros.

\begin{theorem}\label{thm:main}
On $X=\PP^1$ with $\omega=2\omega_{\FS}$, there are smooth potentials $u_0,u_1\in\HH_\omega$ with the following properties. Their weak geodesic $U$ is smooth on a neighborhood of
\[
 (V_-\cup V_+)\times\{\tau:1/4\le\log|\tau|\le3/4\},
\]
where $V_-$ and $V_+$ are neighborhoods of $0$ and $\infty$, respectively. The total form $\Omega_U$ has rank one throughout this neighborhood. For every $t\in[1/4,3/4]$,
\[
 (\omega+\ddc_Xu_t)(0)=(\omega+\ddc_Xu_t)(\infty)=0.
\]
In the local coordinates $\xi=z$ at $0$ and $\xi=1/z$ at $\infty$, the fiber form has the expansion
\begin{equation}\label{eq:expansion}
 \omega+\ddc_Xu_t
 =\left(\frac{16e^{2t}}{(1+e^{2t})^2}|\xi|^2
       +O(|\xi|^4)\right)\ii\,d\xi\wedge d\bar\xi,
\end{equation}
uniformly for $t\in[1/4,3/4]$.
\end{theorem}

Proposition~\ref{prop:explicit} below gives the endpoint potentials and a larger explicit region on which the geodesic agrees with a smooth model. In particular, the zeros in Theorem~\ref{thm:main} persist throughout an open time interval containing $[1/4,3/4]$.

For a compact connected Riemann surface $(C,\omega_C)$, Theorem~\ref{thm:curves} gives endpoints in the prescribed class for which each ramification point of the chosen cover is an isolated fiberwise zero during a common time interval. Its order is $2m-2$ at a point of local degree $m$. By choosing the cover one can realize any prescribed positive even order (Corollary~\ref{cor:curve-prescribed-order}). On a product $C\times Y$ with the product K\"ahler form, the weak geodesic pulls back from $C$. During the same interior time interval, the spatial rank is $n-1$ along each $\{p\}\times Y$ with $p$ a ramification point, while the total rank is $n$ nearby, where $n=\dim_{\mathbb C}(C\times Y)$ (Corollary~\ref{cor:higher-dimensional-products}).

AIMPL Problem~2.5, attributed to Tosatti \cite{aim}, concerns the nondegenerate problems with the same fixed endpoints:
\begin{equation}\label{eq:aim-problem}
 \begin{gathered}
 (\theta+\ddc U_\delta)^{n+1}=\delta\chi^{n+1},
 \qquad\theta+\ddc U_\delta>0,\\
 U_\delta|_{|\tau|=e^j}=u_j\quad(j=0,1),\qquad
 \theta=\pr_X^*\omega,\quad
 \chi=\theta+\ii\,d\tau\wedge d\bar\tau.
 \end{gathered}
\end{equation}
It asks for a counterexample to the assertion that the fiber restrictions satisfy
\[
 \omega+\ddc_XU_\delta(\cdot,\tau)\ge\varepsilon\omega
 \qquad(\tau\in A)
\]
for all sufficiently small $\delta>0$, with $\varepsilon>0$ independent of $\delta$. This is stated as a maximal rank problem with the lower bound imposed on the fiber restrictions. For the endpoints of Theorem~\ref{thm:main}, we obtain the following more localized failure.

\begin{corollary}\label{cor:main-approximation}
For all sufficiently small $\delta>0$, problem~\eqref{eq:aim-problem} has a unique smooth solution for the endpoint pair of Theorem~\ref{thm:main}, and $U_\delta\to U$ uniformly. Fix $t\in[1/4,3/4]$ and a closed coordinate disk $K$ centered at $0$ or $\infty$, sufficiently small that it lies in the corresponding neighborhood of smoothness. If
\[
 \omega+\ddc_XU_\delta(\cdot,e^t)=\lambda_{\delta,t}\omega,
\]
then $\min_K\lambda_{\delta,t}\to0$ as $\delta\downarrow0$. Thus these fixed endpoints supply the counterexample requested in AIMPL Problem~2.5.
\end{corollary}

The proof is given in Section~\ref{sec:consequences}. A fiberwise positive lower bound uniform in $\delta$ would pass to the limit in the sense of currents, contradicting the vanishing of the smooth limit form at the center of the disk. The same argument applies on the open zero set of the one-dimensional examples in \cite[Theorem~A]{xiao-yao}.

\subsection*{Construction and relation to earlier work}
Our starting point is a special case of the smooth log-sum geodesic lines of Darvas--McCleerey \cite[Example~1.4(ii) and Corollary~5.2(ii)]{darvas-mccleerey}. A branched-cover description makes its fiberwise zeros explicit. These zeros occur at every finite time, including at the model endpoints. The main task is therefore to obtain strictly K\"ahler endpoints while preserving the interior degeneration.

We add nonnegative corrections supported near the two ramification points. The endpoint estimates show that the corrected forms are strictly positive. The preservation mechanism is independent of the particular cutoff and numerical constants: Theorem~\ref{thm:preservation} allows different corrections at the two ends, provided they are nonnegative, have the prescribed support, and produce K\"ahler endpoints. Nonnegativity keeps the model below the new weak geodesic. On suitable compact bordered holomorphic curves, the model total form vanishes and the endpoint corrections disappear at the boundary. The maximum principle gives the reverse inequality on these curves and hence exact agreement in the interior region. Locality of the endpoint changes alone would not imply this agreement.

The example complements two relevant constructions. Hu proves positive lower bounds under a smallness assumption on the $C^2$ norm of the boundary data \cite[Theorems~1.1--1.2]{hu}. Ross--Witt Nystr\"om construct complete degeneration for a Dirichlet problem over a disk with boundary data depending on the boundary angle \cite[Theorems~1.1 and~2.1]{ross-maximal}. Here the data are two fixed K\"ahler potentials on the two boundary circles of an annulus. Section~\ref{sec:discussion} explains these distinctions and the obstruction to obtaining the required endpoints by directly truncating the earlier Hele--Shaw ray.

Sections~\ref{sec:model}--\ref{sec:preservation} give the model, endpoint corrections and comparison argument on $\PP^1$. Section~\ref{sec:curves} extends the construction to arbitrary compact connected Riemann surfaces. Section~\ref{sec:consequences} treats the nondegenerate approximations and products.

\section{The semipositive branched model}\label{sec:model}

Normalize $\omega_{\FS}=\ddc\log(1+|z|^2)$ in an affine chart. On $Y=\PP^1_\zeta$ with $\omega_Y=\omega_{\FS}$, consider
\begin{equation}\label{eq:target}
 v_t(\zeta)=\log\frac{1+e^{2t}|\zeta|^2}{1+|\zeta|^2}.
\end{equation}
The total form of $v_{\log|\tau|}$ is the pullback of $\omega_Y$ by $(\zeta,\tau)\mapsto\tau\zeta$. It is semipositive and has square zero, while its fiber restrictions are strictly positive. Thus $v_t$ is a smooth geodesic.

On $X=\PP^1_z$, set
\begin{equation}\label{eq:model-data}
 \omega=2\omega_{\FS},\qquad
 F(z)=\frac{z^2-1}{z^2+1},\qquad
 \wh u_t(z)=\log\frac{|z^2+1|^2+e^{2t}|z^2-1|^2}{(1+|z|^2)^2}.
\end{equation}
In \cite[Example~1.4(ii) and Corollary~5.2(ii)]{darvas-mccleerey}, take the two sections $Z_1^2+Z_0^2$ and $Z_1^2-Z_0^2$ of $\mathcal O(2)$. They have no common zero. Their log-sum potential at time $s=2t$ is precisely $\wh u_t$ in~\eqref{eq:model-data}. The map $F$ is the square map followed by a M\"obius transformation. Its ramification points are $0$ and $\infty$, and their branch values are $-1$ and $1$.

The degree-two pullback $F^*\omega_Y$ is cohomologous to $\omega$. More explicitly, the globally smooth function
\[
 \rho(z)=2\log(1+|z|^2)-\log\bigl(2(1+|z|^4)\bigr)
\]
satisfies
\begin{equation}\label{eq:background}
 \omega=F^*\omega_Y+\ddc\rho,
 \qquad \wh u_t=v_t\circ F-\rho.
\end{equation}
The identities hold away from the poles of $F$ by direct calculation and extend across them by smoothness.

Let $\wh U(z,\tau)=\wh u_{\log|\tau|}(z)$ and
\[
 H:X\times\mathbb C^*\longrightarrow\PP^1,
 \qquad H(z,\tau)=\tau F(z).
\]
In homogeneous coordinates $z=Z_1/Z_0$, these are
\[
 \begin{split}
 H([Z_0:Z_1],\tau)
   &=[Z_1^2+Z_0^2:\tau(Z_1^2-Z_0^2)],\\
 \wh U([Z_0:Z_1],\tau)
   &=\log\frac{|Z_1^2+Z_0^2|^2+|\tau|^2|Z_1^2-Z_0^2|^2}
                    {(|Z_0|^2+|Z_1|^2)^2}.
 \end{split}
\]
The two entries defining $H$ have no common zero. In particular, $\wh U$ is globally smooth for $\tau\ne0$, and
\begin{equation}\label{eq:pullback}
 \wh\Omega:=\pr_X^*\omega+\ddc\wh U
       =H^*\omega_Y\ge0,\qquad \wh\Omega^2=0.
\end{equation}
Consequently, $\wh U$ is an admissible smooth subgeodesic for any endpoint data at least $\wh u_0,\wh u_1$.

Since $F'(z)=4z/(z^2+1)^2$, the fiber form is
\begin{equation}\label{eq:coefficient}
 \omega+\ddc_X\wh u_t
 =g_t(z)\ii\,dz\wedge d\bar z,
 \qquad
 g_t(z)=\frac{16e^{2t}|z|^2}
 {\bigl(|z^2+1|^2+e^{2t}|z^2-1|^2\bigr)^2}.
\end{equation}
The denominator is nonzero everywhere in this chart. The identities $F(1/z)=-F(z)$ and $\rho(1/z)=\rho(z)$ give the same formulas in the coordinate $w=1/z$. Thus the fiber form vanishes precisely at $0$ and $\infty$, and its expansion there is \eqref{eq:expansion}. In particular, no finite segment of the model has strictly K\"ahler endpoints.

\begin{remark}[Total rank and fiber restriction]\label{rem:rank}
The map $H$ is a submersion. Away from the ramification points this follows from $dF\ne0$; at $0$ and $\infty$, the derivative in the $\tau$ direction is nonzero because $F$ equals $-1$ and $1$. Hence $\wh\Omega$ has complex rank one everywhere, the maximal rank allowed by $\wh\Omega^2=0$. At the ramification points its kernel is the spatial tangent line, so its fiber restriction vanishes. On any open set where the corrected geodesic agrees with the model, their total forms and hence their fiber restrictions agree as well.
\end{remark}

\begin{remark}[The standard Hele--Shaw ray]\label{rem:hele-shaw}
The target geodesic \eqref{eq:target} also comes from the standard Hele--Shaw ray with zero initial potential on $Y$:
\[
 G(\zeta,\eta)=\log\frac{|\eta|^2+|\zeta|^2}{1+|\zeta|^2},
 \qquad 0<|\eta|\le1;
\]
see \cite[proof of Lemma~2.5 and Theorem~2.7]{ross-discs}. Indeed,
\[
 G(\zeta,e^{-t})+2t=v_t(\zeta)\qquad(t\ge0).
\]
The complex parameters are related by $\eta=\tau^{-1}$. The added term $2\log|\tau|$ is pluriharmonic on the punctured parameter domain and hence does not change the total form. This ray has a prescribed logarithmic singularity at $(\zeta,\eta)=(0,0)$; the formula is used here only on the punctured disk. Equations~\eqref{eq:background} recover the same branched model. The endpoint correction and interior preservation are the subsequent steps of our construction.
\end{remark}

\section{Repairing the endpoints}\label{sec:endpoints}

We seek corrections $b_0,b_1\ge0$ supported over a narrow modulus band around the branch values of $F$. Theorem~\ref{thm:preservation} will show that every such pair giving K\"ahler endpoints preserves the model on a fixed interior region. This section supplies one explicit choice, with $b_0=b_1=b$.

The local scale is dictated by the quadratic zeros in \eqref{eq:coefficient}. In a coordinate $\xi$ centered at either ramification point, put
\[
 b_R(\xi)=\kappa R^4h(|\xi|^2/R^2)
\]
and extend by zero outside $\{|\xi|<R\}$. Two spatial derivatives give
\[
 (b_R)_{\xi\bar\xi}
   =\kappa R^2\bigl(h'(s)+sh''(s)\bigr),
 \qquad s=|\xi|^2/R^2.
\]
The factor $R^2$ matches the size of the model coefficient where $|\xi|$ is comparable to $R$. We choose $h$ so that its complex Hessian is positive near the center and bounded on the remaining annulus; the model positivity on that annulus then absorbs the negative part of a sufficiently small correction.

\begin{lemma}\label{lem:cutoff}
The function
\begin{equation}\label{eq:cutoff}
 h(s)=
 \begin{cases}
 s\exp\bigl(-s/(1-s)\bigr),&0\le s<1,\\
 0,&s\ge1
 \end{cases}
\end{equation}
is smooth and nonnegative on $[0,\infty)$, and all its derivatives vanish at $1$. For $k=h'+sh''$,
\begin{equation}\label{eq:k}
 k(s)=e^{-s/(1-s)}
 \frac{s^4-5s^3+11s^2-7s+1}{(1-s)^4}
 \qquad(0\le s<1).
\end{equation}
Moreover, $k(0)=1$, $k(s)>0$ for $0\le s\le1/10$, and $|k(s)|<64$ for all $s\ge0$.
\end{lemma}

\begin{proof}
The exponential factor gives smoothness and vanishing of all derivatives at $1$, and differentiation gives \eqref{eq:k}. On $[0,1/10]$, its numerator is at least
\[
 1-7s-5s^3\ge1-7/10-5/1000>0.
\]
For the global bound, set $x=(1-s)^{-1}\ge1$. Formula~\eqref{eq:k} gives
\[
 |k(s)|\le e^{1-x}(1+3x^2+2x^3+x^4).
\]
The maxima of $e^{1-x}x^j$ for $j=2,3,4$ on $[1,\infty)$ are $4/e<2$, $27/e^2<7$, and $256/e^3<32$. Thus $|k|<1+6+14+32<64$ for $s<1$; for $s\ge1$, it is zero.
\end{proof}

Fix
\begin{equation}\label{eq:band}
 a=\sqrt{19/20},\qquad d=\sqrt{21/20},\qquad
 W=\{\zeta\in\mathbb C:a<|\zeta|<d\}.
\end{equation}
The branch values $-1$ and $1$ lie in this same band. Taking $R=1/10$ and $\kappa=10^{-4}$ in the two disjoint ramification charts gives
\begin{equation}\label{eq:bump}
 b(z)=10^{-8}\bigl[h(100|z|^2)+h(100|z|^{-2})\bigr],
 \qquad b(0)=b(\infty)=0.
\end{equation}
Near $0$ and $\infty$, the two summands are interpreted in the coordinates $z$ and $w=1/z$, respectively; each is extended by zero outside its coordinate disk.

\begin{proposition}\label{prop:endpoints}
The function $b$ is smooth and nonnegative, and
\begin{equation}\label{eq:support}
 \supp b\Subset F^{-1}(W),\qquad \norm b_{C^0}\le10^{-8}.
\end{equation}
The endpoint potentials
\begin{equation}\label{eq:endpoints}
 u_0=\wh u_0+b,\qquad u_1=\wh u_1+b
\end{equation}
belong to $\HH_\omega$.
\end{proposition}

\begin{proof}
Lemma~\ref{lem:cutoff} gives smoothness at the cutoff radii. The supports of the two terms in \eqref{eq:bump} are disjoint, and $0\le h\le1$, giving the $C^0$ bound. If $|z|\le1/10$, then
\[
 \frac{99}{101}\le|F(z)|\le\frac{101}{99},\qquad
 \frac{19}{20}<\left(\frac{99}{101}\right)^2
 \le|F(z)|^2\le\left(\frac{101}{99}\right)^2<\frac{21}{20}.
\]
The same estimate holds in the $w$ chart, so the support is compactly contained in $F^{-1}(W)$.

For endpoint positivity near $0$, put $r=|z|^2\le1/100$ and $q=e^{2j}$, $j\in\{0,1\}$. Then
\[
 |z^2+1|^2+q|z^2-1|^2\le(1+q)(1+r)^2.
\]
Since $1\le q\le e^2<8$ and $q/(1+q)^2$ is nonincreasing for $q\ge1$, \eqref{eq:coefficient} yields
\begin{equation}\label{eq:model-bound}
 g_j(z)\ge\frac{128}{81}\left(\frac{100}{101}\right)^4r\ge r.
\end{equation}
The final numerical factor is greater than one. In this chart,
\[
 b_{z\bar z}=10^{-6}k(100r).
\]
This is positive for $100r\le1/10$, including at $r=0$. In the remainder of the support, $r\ge10^{-3}$, and therefore
\[
 g_j+b_{z\bar z}\ge10^{-3}-64\cdot10^{-6}>0.
\]
The identical argument in the $w$ chart treats $\infty$. Outside the two disks the correction vanishes and the model fiber forms are positive.
\end{proof}

\begin{remark}[Size of the correction]\label{rem:small}
The correction can be arbitrarily small in $C^2$. For sufficiently small $R$, the model satisfies $g_j\ge c|\xi|^2$ near each ramification point. A sufficiently small fixed $\kappa>0$ then makes the corrected endpoints positive by the core and annulus estimates above. In a fixed smooth background norm,
\[
 \norm{b_R}_{C^0}=O(R^4),\qquad \norm{b_R}_{C^2}=O(R^2).
\]
This controls the change from the semipositive model. The size of the full endpoint data in small-data positivity results is discussed in Section~\ref{sec:discussion}.
\end{remark}

\section{Exact preservation in the interior}\label{sec:preservation}

The model's null curves explain why the endpoint repair leaves an interior region unchanged. A curve on which $H=\tau F=c$ is constant is null for $\wh\Omega$ by \eqref{eq:pullback}. Its intersections with $|\tau|=1$ and $|\tau|=e$ occur at the spatial levels $|F|=|c|$ and $|F|=|c|/e$. Both levels avoid the support band $W$ if
\[
 |c|>d,\qquad |c|/e<a.
\]
This leads to the open region
\begin{equation}\label{eq:region}
 \RR=\{(z,\tau)\in X\times A:d<|\tau F(z)|<ea\}.
\end{equation}
Only finite nonzero affine values of $F$ occur there. Ross and Witt Nystr\"om use a related comparison along harmonic discs in the proof of their Theorem~4.2 \cite[p.~332]{ross-discs}. Here the comparison takes place on compact bordered Riemann surfaces whose projections to the time annulus may be branched.

\begin{lemma}[Comparison on a null curve]\label{lem:null-comparison}
Let $(X,\omega)$ be a compact K\"ahler manifold and $U$ the continuous weak geodesic between two smooth K\"ahler endpoints. Let $V$ be a smooth real function on a neighborhood of $X\times\overline A$ such that $\pr_X^*\omega+\ddc V\ge0$ on $X\times A$ and $V\le U$ on $X\times\partial A$. Suppose $L$ is a compact connected bordered Riemann surface with nonempty boundary and $j$ is holomorphic on a neighborhood of $L$, with
\[
 j(L)\subset X\times\overline A,\qquad
 j(L^\circ)\subset X\times A,\qquad
 j(\partial L)\subset X\times\partial A.
\]
If
\[
 j^*(\pr_X^*\omega+\ddc V)=0\quad\text{on }L^\circ,
 \qquad (U-V)\circ j=0\quad\text{on }\partial L,
\]
then $U\circ j=V\circ j$ on $L$.
\end{lemma}

\begin{proof}
Set $V^*(x,\tau)=\max_{|\gamma|=1}V(x,\gamma\tau)$. This maximum is continuous by compactness of the circle. After adding a local potential for $\omega$ pulled back from $X$, each rotated function is plurisubharmonic, and their continuous upper envelope is plurisubharmonic as well. Thus $V^*$ is a rotation-invariant subgeodesic. The endpoint values of $U$ are constant along each boundary circle for fixed $x$, so $V^*\le U$ there. The envelope characterization \cite[Lemma~2.1]{darvas-survey} gives $U\ge V^*\ge V$ throughout. Set $f=(U-V)\circ j$.
Around each point of $L^\circ$, take a coordinate disk $D$ small enough
that $\pr_X(j(D))$ lies in a chart with a local potential $p$ for $\omega$.
Then $p\circ\pr_X+U$ is continuous and plurisubharmonic near $j(D)$.
By the holomorphic pullback theorem \cite[Ch.~I, Definition~(5.1) and
Theorem~(5.11)]{demailly}, its composition with $j|_D$ is
plurisubharmonic on $D$, hence subharmonic in the disk coordinate.
This applies in particular when $j(D)$ lies on the graph in~\eqref{eq:null-curve}.
The null condition gives
$\ddc((p\circ\pr_X+V)\circ j)=j^*(\pr_X^*\omega+\ddc V)=0$ on $D$,
so this second composition is harmonic. Their difference $f$ is
therefore subharmonic on $L^\circ$. It is continuous on $L$, nonnegative, and zero on $\partial L$. The maximum principle gives $f=0$.
\end{proof}

We describe the bordered curves for any nonconstant holomorphic map
$F:C\to\PP^1$ from a compact connected Riemann surface. Fix
$c\in\mathbb C^*$ such that the circles $|w|=|c|/e$ and $|w|=|c|$
avoid all branch values, and set
\begin{equation}\label{eq:null-curve}
 K_c=\{z\in C: |c|/e\le |F(z)|\le |c|\},\qquad
 j_c(z)=(z,c/F(z)).
\end{equation}
In applying Lemma~\ref{lem:null-comparison}, we take the map $j$ on the
chosen connected component $L$ of $K_c$ to be the restriction $j_c|_L$.
The two boundary levels are regular levels of $\log|F|$ away from the
zeros and poles of $F$, so each
component $L$ of $K_c$ is a compact Riemann surface with smooth
boundary. A component with empty boundary would be open and closed in
$C$, whereas $F$ is surjective and $K_c$ excludes its zeros and poles.
The graph $j_c$ is holomorphic on a neighborhood of $L$ and remains
an immersion at ramification points of $F$. It maps $L^\circ$ into
$C\times A$ and the levels $|F|=|c|,|c|/e$ to the circles
$|\tau|=1,e$, respectively. Along this graph $\tau F(z)=c$.

We now apply this construction to $X=\PP^1$ and $F$ in
\eqref{eq:model-data}, with the band $W$ from \eqref{eq:band}.

\begin{theorem}[Preservation under endpoint corrections]\label{thm:preservation}
Let $b_0,b_1\in C^\infty(X)$ be nonnegative and satisfy
\[
 \supp b_j\Subset F^{-1}(W),\qquad
 \wh u_j+b_j\in\HH_\omega\quad(j=0,1).
\]
Then the weak geodesic joining $\wh u_0+b_0$ and $\wh u_1+b_1$ agrees with $\wh U$ on $\RR$. The conclusion requires no additional norm bound on $b_0$ or $b_1$.
\end{theorem}

\begin{proof}
Write $U$ for this geodesic. Nonnegativity of $b_0,b_1$ and
\eqref{eq:pullback} make $\wh U$ admissible for the raised endpoint data.
Fix $(z_0,\tau_0)\in\RR$ and put $c=\tau_0F(z_0)$. Since
\[
 |c|/e<a<1<d<|c|,
\]
both boundary circles avoid the branch values $-1,1$. Take the component
$L$ of $K_c$ in \eqref{eq:null-curve} containing $z_0$. The identity
$H\circ j_c=c$ gives $j_c^*\wh\Omega=0$. Both boundary levels lie outside
$F^{-1}(W)$, so the corrections vanish there and
$(U-\wh U)\circ j_c=0$ on $\partial L$. Lemma~\ref{lem:null-comparison}
gives equality at $j_c(z_0)=(z_0,\tau_0)$.
\end{proof}

\begin{proposition}[An explicit endpoint pair]\label{prop:explicit}
For the function $b$ in \eqref{eq:bump}, the endpoints \eqref{eq:endpoints} lie in $\HH_\omega$, and their weak geodesic satisfies
\begin{equation}\label{eq:agreement}
 U=\wh U\qquad\text{on }\RR.
\end{equation}
In particular, it is smooth there, and the fiber form vanishes at $0$ and $\infty$ for every
\begin{equation}\label{eq:time-interval}
 t\in I:=(\log d,1+\log a).
\end{equation}
\end{proposition}

\begin{proof}
Proposition~\ref{prop:endpoints} verifies the hypotheses of Theorem~\ref{thm:preservation}, giving \eqref{eq:agreement}. At the ramification points, $|F|=1$, so membership in $\RR$ is equivalent to $t\in I$. Formula~\eqref{eq:coefficient} gives the vanishing and expansion~\eqref{eq:expansion}, uniformly on compact subintervals of $I$.
\end{proof}

\begin{proof}[Proof of Theorem~\ref{thm:main}]
The inequalities
\[
 e^{1/4}a>d,\qquad e^{3/4}d<ea
\]
are both equivalent to $e^{1/2}>21/19$. Therefore
\[
 F^{-1}(W)\times\{\tau:1/4\le\log|\tau|\le3/4\}\subset\RR.
\]
Choose $V_-$ and $V_+$ around $0$ and $\infty$ with closures in $F^{-1}(W)$. Since $\RR$ is open, Proposition~\ref{prop:explicit} gives smoothness on a neighborhood of their product with the closed time interval. There $\Omega_U=H^*\omega_{\FS}$ has rank one by Remark~\ref{rem:rank}. The expansion follows from Proposition~\ref{prop:explicit}.
\end{proof}

\begin{remark}
The comparison maps $j_c$ are graphs over the spatial variable. Their projections to the time annulus are branched at the ramification points of $F$, precisely where the model kernel becomes spatial. The comparison uses the bordered curves themselves and does not require a foliation theorem for the unknown geodesic.
\end{remark}

\section{Extension to compact Riemann surfaces}\label{sec:curves}

The same construction works on every compact connected Riemann surface,
with its background K\"ahler form prescribed in advance. The two additional
steps are to place all branch values in a common modulus band and to match
the area of the pullback form to that of the given background. The null-curve
comparison in Lemma~\ref{lem:null-comparison} then applies without change.

\begin{theorem}\label{thm:curves}
Let $C$ be a compact connected Riemann surface and let $\omega_C$ be any
smooth K\"ahler form on $C$. There are smooth endpoint potentials
$\varphi_0,\varphi_1\in\HH_{\omega_C}$, a nonconstant holomorphic map
$F:C\to\PP^1$, and a nonempty finite set $B=\operatorname{Crit}(F)$ with
the following properties. The weak geodesic $U$ joining $\varphi_0$ and $\varphi_1$ is smooth on an open set
$\RR_C\subset C\times A$ containing
\[
 B\times\{\tau:3/8<\log|\tau|<5/8\}.
\]
On $\RR_C$, the total form
\[
 \Omega_U^C:=\pr_C^*\omega_C+\ddc U
\]
has complex rank one. If $u_t=U(\cdot,e^t)$, then
\begin{equation}\label{eq:curve-fiber-zero}
 (\omega_C+\ddc_Cu_t)(p)=0
 \qquad(p\in B,\ 3/8<t<5/8).
\end{equation}
At a point $p\in B$ of ramification index $m_p\ge2$, this fiber form has a
zero of order exactly $2m_p-2$ in a centered local coordinate. The leading
coefficient is positive, and the local expansion is uniform for $t$ in
compact subintervals of $(3/8,5/8)$.
\end{theorem}

\subsection{Normalizing the branch values}

\begin{lemma}\label{lem:curve-normalization}
There is a nonconstant holomorphic map $F:C\to\PP^1$ of degree $d\ge2$
such that $0$ and $\infty$ are regular values and
\begin{equation}\label{eq:curve-critical-band}
 e^{-1/8}<|F(p)|<e^{1/8}
 \qquad\text{for every }p\in B:=\operatorname{Crit}(F).
\end{equation}
The set $B$ is finite and nonempty.
\end{lemma}

\begin{proof}
Let $g=g(C)$ and fix $q\in C$. By Riemann--Roch, the space of meromorphic
functions with no poles away from $q$ and pole order at most $g+1$ at
$q$ has dimension at least $(g+1)+1-g=2$. Since the constants span only
one dimension, this space contains a nonconstant function $f$. Setting
$f=\infty$ at a pole makes $f$ a holomorphic map $C\to\PP^1$, since
in the target coordinate $1/w$ the function $1/f$ is holomorphic and
vanishes at the pole. If $f$ has degree one, postcompose
it with $w\mapsto w^2$. Write $f$ for the resulting map and $d\ge2$
for its degree. The Riemann--Hurwitz formula gives
\[
 \sum_{p\in C}(e_p-1)=2g(C)-2+2d>0,
\]
where $e_p$ is the local degree of $f$ at $p$. Hence its critical set is
nonempty; it is finite because $f$ is nonconstant and $C$ is compact.

Let $S$ be the finite set of critical values of $f$, and fix
$a\in\mathbb C\setminus S$. Choose $0<\eta<\min\{1-e^{-1/8},e^{1/8}-1\}$,
so that $(1-\eta,1+\eta)\subset(e^{-1/8},e^{1/8})$.
Since $S$ is finite and $a\notin S$, we can choose $r>0$ so that
$|s-a|\ge 2r$ for every finite $s\in S$. Take any $b\in\mathbb C$
with $0<|b-a|<r\eta$. As $\eta<1$, the disk $|w-a|<r$ contains $b$
and no critical values, so $b$ is a regular value distinct from $a$.

For $T_{a,b}(w)=(w-a)/(w-b)$ and every finite $s\in S$, we now have
\[
 |s-b|\ge |s-a|-|b-a|>r,\qquad
 |T_{a,b}(s)-1|=\frac{|b-a|}{|s-b|}<\eta .
\]
The reverse triangle inequality gives
$\bigl||T_{a,b}(s)|-1\bigr|<\eta$. If $\infty\in S$, its image is
$T_{a,b}(\infty)=1$. Thus the same choice of $b$ puts every element of
$T_{a,b}(S)$ in the desired modulus band. Since $T_{a,b}$ is a
M\"obius transformation, $F=T_{a,b}\circ f$ has the same ramification
points and local degrees as $f$, with critical values $T_{a,b}(S)$.
This proves \eqref{eq:curve-critical-band}. Moreover,
$F^{-1}(0)=f^{-1}(a)$ and $F^{-1}(\infty)=f^{-1}(b)$, so both $0$ and
$\infty$ are regular values of $F$.
The M\"obius-normalization step applies to any given map of degree
at least two and preserves all its local degrees.
\end{proof}

\subsection{The model in the prescribed background class}

Fix $F$ as in Lemma~\ref{lem:curve-normalization}. With the normalization
$\omega_{\FS}=\ddc\log(1+|w|^2)$ used above, put
\begin{equation}\label{eq:curve-area-normalization}
 \lambda=
 \frac{\int_C\omega_C}{d\int_{\PP^1}\omega_{\FS}}>0.
\end{equation}
The degree formula gives $\int_C F^*\omega_{\FS}=d\int_{\PP^1}\omega_{\FS}$.
Thus the smooth real two-form $\sigma=\lambda F^*\omega_{\FS}-\omega_C$
has integral zero. On the compact connected oriented surface $C$,
integration identifies $H^2_{\mathrm{dR}}(C;\mathbb R)$ with $\mathbb R$,
so $\sigma$ is exact. To obtain a global potential, write $\sigma=h\omega_C$.
Then $\int_C h\omega_C=0$. With the normalization
$\ddc u=(\Delta_{\omega_C}u)\omega_C$, the Poisson equation
$\Delta_{\omega_C}\rho_C=h$ therefore has a smooth real solution, unique
up to an additive constant. It satisfies
\begin{equation}\label{eq:curve-background}
 \omega_C+\ddc\rho_C=\lambda F^*\omega_{\FS}.
\end{equation}
Thus the scaled pullback and $\omega_C$ represent the same cohomology
class, and \eqref{eq:curve-background} exhibits $\rho_C$ as a global
potential for their difference. The pullback vanishes at the ramification points,
where $\omega_C$ is strictly positive.

Define
\begin{equation}\label{eq:curve-model}
 \begin{split}
 V(z,\tau)&=\rho_C(z)+\lambda\log
       \frac{1+|\tau F(z)|^2}{1+|F(z)|^2},\\
 H_C(z,\tau)&=\tau F(z),
 \qquad (z,\tau)\in C\times\mathbb C^*.
 \end{split}
\end{equation}
The map $H_C:C\times\mathbb C^*\to\PP^1$ is holomorphic. The function
$V$ is globally smooth: near a pole of $F$, the coordinate $q=1/F$ gives
\[
 \log\frac{1+|\tau F|^2}{1+|F|^2}
 =\log|\tau|^2+
   \log\frac{1+|q/\tau|^2}{1+|q|^2}.
\]
In particular, $V$ extends smoothly to a neighborhood of
$C\times\overline A$. If $\theta_C=\pr_C^*\omega_C$, then
\eqref{eq:curve-background} gives
\begin{equation}\label{eq:curve-pullback}
 \Omega_V^C:=\theta_C+\ddc V
   =\lambda H_C^*\omega_{\FS}\ge0,
 \qquad (\Omega_V^C)^2=0.
\end{equation}
This follows in the affine target chart by differentiation and extends
across the poles by smoothness. The square vanishes because the target has
complex dimension one. Thus $V$ is a smooth $S^1$-invariant subgeodesic.

Write $v_j=V(\cdot,e^j)$ for $j=0,1$. Their fiber forms are
\begin{equation}\label{eq:curve-model-endpoints}
 \eta_j:=\omega_C+\ddc v_j
       =\lambda(e^jF)^*\omega_{\FS}.
\end{equation}
They vanish exactly at $B$ and are strictly positive on $C\setminus B$.

\subsection{Nonnegative endpoint corrections}

Set
\begin{equation}\label{eq:curve-support-band}
 W_C=\{w\in\mathbb C:e^{-1/4}<|w|<e^{1/4}\}.
\end{equation}
By \eqref{eq:curve-critical-band}, $B\subset F^{-1}(W_C)$.

\begin{lemma}\label{lem:curve-repair}
There is a smooth real function $\beta$ on $C$ such that
\[
 \beta\ge0,\qquad \supp\beta\Subset F^{-1}(W_C),\qquad
 \varphi_j:=v_j+\beta\in\HH_{\omega_C}\quad(j=0,1).
\]
\end{lemma}

\begin{proof}
Choose pairwise disjoint coordinate disks around the finitely many points
of $B$, with closures in $F^{-1}(W_C)$. At $p\in B$, take a centered
coordinate $\xi$ and a radius $r_p>0$ such that the closed disk
$|\xi|\le2r_p$ lies in its coordinate neighborhood. Choose a smooth
function $0\le\chi_p\le1$ equal to one on $|\xi|\le r_p$ and compactly
supported in $|\xi|<2r_p$. For a constant $\epsilon_p>0$, put
\[
 \beta_p=\epsilon_p\chi_p|\xi|^2,
\]
and extend it by zero to $C$.

Write $\eta_j=\ii g_j\,d\xi\wedge d\bar\xi$. On the inner disk,
$g_j\ge0$ and $(\beta_p)_{\xi\bar\xi}=\epsilon_p>0$. On the closed
transition annulus the only critical point in the coordinate disk has
been excluded, so
\[
 \mu_p:=\min_{j=0,1}\ \min_{r_p\le|\xi|\le2r_p}g_j>0.
\]
Choose
\[
 K_p=\max\left\{1,
 \sup_{|\xi|\le2r_p}
  \left|\partial_\xi\partial_{\bar\xi}(\chi_p|\xi|^2)\right|
 \right\},
 \qquad 0<\epsilon_p<\frac{\mu_p}{2K_p}.
\]
Then $g_j+(\beta_p)_{\xi\bar\xi}\ge\mu_p/2>0$ on the transition
annulus, for both $j=0,1$. The zero extension is smooth because $\chi_p$
has compact support in the coordinate disk. Set
$\beta=\sum_{p\in B}\beta_p$. The supports are disjoint, and outside
these disks the model forms are already strictly positive. The inner-disk
and transition-annulus estimates therefore prove
$\eta_j+\ddc\beta>0$ on all of $C$. Nonnegativity and the asserted
support property follow directly from the construction.
\end{proof}

\subsection{Preservation on compact bordered null curves}

The corresponding interior region is
\begin{equation}\label{eq:curve-region}
 \RR_C=\{(z,\tau)\in C\times A:
                e^{1/4}<|\tau F(z)|<e^{3/4}\}.
\end{equation}
Here and below this inequality selects only finite nonzero values of $F$;
the region is open in $C\times A$.

\begin{proposition}\label{prop:curve-preservation}
Let $\beta_0,\beta_1\in C^\infty(C,\mathbb R)$ satisfy
\[
 \beta_j\ge0,\qquad \supp\beta_j\Subset F^{-1}(W_C),\qquad
 v_j+\beta_j\in\HH_{\omega_C}\quad(j=0,1).
\]
The weak geodesic joining $v_0+\beta_0$ and $v_1+\beta_1$ agrees with
$V$ on $\RR_C$.
\end{proposition}

\begin{proof}
Write $U$ for this geodesic. By \eqref{eq:curve-pullback} and the
nonnegativity of $\beta_0,\beta_1$, the model $V$ is admissible for
the raised endpoints. Fix $(z_*,\tau_*)\in\RR_C$ and put
$c=\tau_*F(z_*)$. By \eqref{eq:curve-region},
\begin{equation}\label{eq:curve-boundary-separation}
 |c|/e<e^{-1/4}<e^{-1/8},
 \qquad |c|>e^{1/4}>e^{1/8}.
\end{equation}
By \eqref{eq:curve-critical-band}, both boundary circles avoid all branch
values of $F$. Apply the construction \eqref{eq:null-curve} to the component
$L$ of $K_c$ containing $z_*$. The graph satisfies $H_C\circ j_c=c$, so
\eqref{eq:curve-pullback} gives $j_c^*\Omega_V^C=0$. By
\eqref{eq:curve-boundary-separation}, both boundary levels are disjoint
from $F^{-1}(W_C)$, so $(U-V)\circ j_c=0$ on $\partial L$.
Lemma~\ref{lem:null-comparison} now gives equality at
$j_c(z_*)=(z_*,\tau_*)$.
\end{proof}

\subsection{Fiberwise vanishing and total rank}

\begin{proof}[Proof of Theorem~\ref{thm:curves}]
Take the common correction $\beta$ from Lemma~\ref{lem:curve-repair} and
let $U$ join $\varphi_j=v_j+\beta$. The weak geodesic existence and
continuity statements recalled in the introduction apply to
$(C,\omega_C)$. Proposition~\ref{prop:curve-preservation} gives $U=V$ on
$\RR_C$. If $p\in B$ and $3/8<t<5/8$, then
\eqref{eq:curve-critical-band} implies
\[
 \frac14<t+\log|F(p)|<\frac34.
\]
Thus $B\times\{\tau:3/8<\log|\tau|<5/8\}\subset\RR_C$.
The smoothness of $V$ proves the asserted smoothness of $U$ there. On
the part of a fixed fiber contained in this region,
\begin{equation}\label{eq:curve-fiber-pullback}
 \omega_C+\ddc_Cu_t=\lambda(e^tF)^*\omega_{\FS}.
\end{equation}
In particular the fiber form vanishes at every $p\in B$.

Fix $p\in B$ and $t\in(3/8,5/8)$, write $m=m_p$, and take a coordinate
$\xi$ centered at $p$.
Since $F(p)$ is finite and nonzero, its local expansion is
\[
 F(\xi)=a_p+\alpha_p\xi^m+O(\xi^{m+1}),
 \qquad a_p\alpha_p\ne0.
\]
Differentiating the Taylor series gives $F'(\xi)=\xi^{m-1}h_p(\xi)$,
where $h_p$ is holomorphic and $h_p(0)=m\alpha_p\ne0$.
On the fixed time slice, write
$\omega_C+\ddc_Cu_t=g_t(\xi)\ii\,d\xi\wedge d\bar\xi$.
By \eqref{eq:curve-fiber-pullback},
\[
 g_t(\xi)=|\xi|^{2m-2}A_t(\xi),\qquad
 A_t(\xi)=\frac{\lambda e^{2t}|h_p(\xi)|^2}
                   {(1+e^{2t}|F(\xi)|^2)^2}.
\]
The factor $A_t$ is smooth and $A_t(0)>0$. In particular,
\begin{equation}\label{eq:curve-vanishing-order}
 g_t(\xi)=
 \frac{\lambda e^{2t}m^2|\alpha_p|^2}
      {(1+e^{2t}|a_p|^2)^2}|\xi|^{2m-2}
       +O(|\xi|^{2m-1}).
\end{equation}
For a compact interval $J\Subset(3/8,5/8)$, the coordinate disk can be
chosen small enough that its product with
$\{\tau:\log|\tau|\in J\}$ lies in $\RR_C$. The expansion is then
uniform for $t\in J$. The factorization shows that for each
$t\in(3/8,5/8)$ the zero at the same point $p$ has order exactly
$2m-2$ in the spatial variable $\xi$. It is isolated since $h_p$
is nonzero near the origin.

Finally, $F$ is finite and nonzero throughout $\RR_C$, so
$\partial_\tau H_C=F\ne0$ there. Hence $H_C$ is a submersion on this
region, and $\Omega_U^C=\lambda H_C^*\omega_{\FS}$ has complex rank one.
\end{proof}

\begin{corollary}[Prescribing an even vanishing order]
\label{cor:curve-prescribed-order}
For every compact connected Riemann surface $C$, every smooth K\"ahler
form $\omega_C$, and every integer $m\ge2$, the endpoint pair in
Theorem~\ref{thm:curves} can be chosen so that, at some fixed point
$p\in C$, the fiber form has a zero of order exactly $2m-2$ for every
$t\in(3/8,5/8)$.
\end{corollary}

\begin{proof}
Take a nonconstant meromorphic function $f:C\to\PP^1$, choose a finite
regular value $a$, and fix $p\in f^{-1}(a)$. Since $df(p)\ne0$,
$\xi=f-a$ is a local coordinate centered at $p$. The meromorphic map
\[
 G=(f-a)^m:C\longrightarrow\PP^1
\]
has degree $m\deg f\ge2$ and satisfies $G(\xi)=\xi^m$ near $p$.
Thus its local degree at $p$ is exactly $m$.

We now apply the M\"obius-normalization step in the proof of
Lemma~\ref{lem:curve-normalization} to this particular map $G$.
Let $S_G$ be its finite set of critical values. Choose a finite regular
value $\alpha$ of $G$, and a distinct finite regular value $\beta$
sufficiently close to $\alpha$, and put
\[
 T(w)=\frac{w-\alpha}{w-\beta},\qquad F=T\circ G.
\]
Since $T$ is a M\"obius transformation, the critical values of $F$ are
exactly $T(S_G)$. For every finite $s\in S_G$,
$|T(s)-1|=|\beta-\alpha|/|s-\beta|$, while $T(\infty)=1$.
As in the proof of the lemma, finiteness of $S_G$ allows one choice of
$\beta$ for which all critical values of $F$ lie in
$e^{-1/8}<|w|<e^{1/8}$. Moreover, $0$ and $\infty$ are regular values
of $F$, since their inverse images under $T$ are the regular values
$\alpha$ and $\beta$ of $G$.

Postcomposition by the M\"obius transformation preserves every local
degree. In particular, this can be checked directly at our chosen point.
Since $G(p)=0$ is a critical value, $\alpha\beta\ne0$, and
\[
 F(\xi)=\frac{\xi^m-\alpha}{\xi^m-\beta}
       =\frac{\alpha}{\beta}
        +\frac{\alpha-\beta}{\beta^2}\xi^m
        +O(\xi^{2m}).
\]
The coefficient of $\xi^m$ is nonzero, so $F$ still has local degree
exactly $m$ at $p$, with $F(p)$ finite and nonzero.

Carry out the background normalization, endpoint correction and interior
preservation construction of Theorem~\ref{thm:curves} using this $F$.
It gives one fixed pair of smooth K\"ahler endpoints whose weak geodesic
agrees with the model near $p$ for every $t\in(3/8,5/8)$.
By~\eqref{eq:curve-vanishing-order}, local degree $m$ gives a fiber zero
of order exactly $2m-2$ at this same point $p$ throughout the interval.
No condition on the local degrees at the other critical points is needed.
\end{proof}

\section{Consequences for positivity and approximation}
\label{sec:consequences}

Let $C$ be a compact connected Riemann surface with a prescribed
K\"ahler form $\omega$, and write $X=C$. The approximation proposition
below applies to any fixed smooth strictly K\"ahler endpoints $u_0,u_1$;
write $U$ for their weak geodesic and $u_t=U(\cdot,e^t)$.
For the degeneration statements, we use either of the two constructions
above. For Proposition~\ref{prop:explicit}, put $J=I$ and
$B=\{0,\infty\}$; for Theorem~\ref{thm:curves}, put
$J=(3/8,5/8)$ and let $B$ be its finite nonempty ramification set,
denoting its endpoints also by $u_0,u_1$.
In both cases, local smoothness and fiberwise zeros exclude a positive
lower bound and force degeneration of the approximating fiber forms.

\begin{corollary}
\label{cor:darvas-obstruction}
For either constructed endpoint pair, every $t\in J$, and every
$\varepsilon>0$,
\[
  \omega+\ddc_Xu_t\not\geq\varepsilon\omega
\]
as currents on $X$. In particular, the strict-positivity question recorded
in Darvas's survey
\cite[Question~2.3]{darvas-survey} has a negative answer.
\end{corollary}

\begin{proof}
Fix $t\in J$ and $p\in B$. By Proposition~\ref{prop:explicit} or
Theorem~\ref{thm:curves}, the fiber form is smooth near $p$ and
vanishes at $p$. An inequality between these smooth forms in the sense of
currents would hold pointwise there, contradicting the positivity of
$\omega$ at $p$.
\end{proof}

\subsection{Nondegenerate approximations with fixed endpoints}
\label{subsec:fixed-endpoint-approximation}

Set
\[
  \theta=\pr_X^*\omega,
  \qquad
  \chi=\theta+\ii\,d\tau\wedge d\bar\tau.
\]
We consider the Dirichlet problems
\begin{equation}
\label{eq:nondegenerate-dirichlet}
  (\theta+\ddc U_\delta)^2=\delta\chi^2,
  \qquad
  \theta+\ddc U_\delta>0,
  \qquad
  U_\delta\big|_{|\tau|=e^j}=u_j\quad(j=0,1).
\end{equation}
These are the complex-dimension-one problems in AIMPL
Problem~2.5~\cite{aim}.
The AIMPL normalization is exactly the one in~\eqref{eq:nondegenerate-dirichlet}.
Since $\chi^2=2\theta\wedge\ii\,d\tau\wedge d\bar\tau$, the parameter
in Hu's equation \cite[(1.14)]{hu} is $\epsilon=2\delta$ in our notation. The strict
positivity in~\eqref{eq:nondegenerate-dirichlet} gives a positive lower
bound for the fiber restrictions of each individual solution. The issue
is whether such a bound can be independent of $\delta$.

\begin{proposition}
\label{prop:approximation-convergence}
For any fixed $u_0,u_1\in\HH_\omega$ on $C$, there are constants
$\alpha>0$ and $C_0<\infty$, depending only on these fixed data, such
that for every $0<\delta\leq\alpha$,
\eqref{eq:nondegenerate-dirichlet} has a unique smooth solution on
$X\times\overline A$. This solution is $S^1$-invariant and satisfies
\begin{equation}
\label{eq:uniform-approximation}
  0\leq U-U_\delta\leq C_0\sqrt\delta
  \qquad\text{on }X\times\overline A.
\end{equation}
\end{proposition}

\begin{proof}
We first construct a strict subsolution without changing the endpoint
values. Write $f=u_1-u_0$ and set
\[
  S(x,\tau)=u_0(x)+tf(x)+L(t^2-t),
  \qquad t=\log|\tau|.
\]
The fiber block of $\theta+\ddc S$ is
\[
  \omega_t^{\mathrm{lin}}
    =(1-t)(\omega+\ddc u_0)+t(\omega+\ddc u_1),
\]
which has a positive lower bound on $X\times[0,1]$. In a local time
coordinate $\zeta=\log\tau$, write
$\omega_t^{\mathrm{lin}}=\ii g_t^{\mathrm{lin}}\,dz\wedge d\bar z$.
The Hermitian matrix of the total form is
\[
  \begin{pmatrix}
    g_t^{\mathrm{lin}}&f_z/2\\
    f_{\bar z}/2&L/2
  \end{pmatrix}.
\]
Thus its Schur complement is
\[
  \frac L2-\frac14|\partial f|_{\omega_t^{\mathrm{lin}}}^2.
\]
Choose, once and for all,
\[
  L>\frac12\sup_{X\times[0,1]}
       |\partial f|_{\omega_t^{\mathrm{lin}}}^2.
\]
Then $\eta:=\theta+\ddc S$ is a smooth K\"ahler form on the closed product,
and $S$ has boundary values $u_0,u_1$. This is an explicit instance of the
strict extension construction in~\cite[Proposition~7.10]{boucksom}.
Its auxiliary-factor hypothesis holds because $q(\tau)=t^2-t$ vanishes
on both boundary circles and satisfies
$\ddc q=\ii\,d\tau\wedge d\bar\tau/(2|\tau|^2)>0$ on $\overline A$.

Choose $\alpha>0$ such that $\eta^2\geq\alpha\chi^2$ everywhere on
$X\times\overline A$. For $0<\delta\leq\alpha$, the substitution
$\psi_\delta=U_\delta-S$ gives
\[
  (\eta+\ddc\psi_\delta)^2=\delta\chi^2,
  \qquad \eta+\ddc\psi_\delta>0,
  \qquad \psi_\delta|_{\partial(X\times\overline A)}=0.
\]
The zero function is a smooth subsolution, since
$\eta^2\geq\delta\chi^2$. Thus the background is K\"ahler, the
right-hand side is a smooth positive volume form, and a smooth subsolution
with the prescribed boundary values is available. The nondegenerate
Dirichlet theorem~\cite[Theorem~A]{boucksom} gives a unique smooth solution. Rotations of $\tau$ preserve the equation and
boundary values, so uniqueness also gives $S^1$-invariance.

The envelope characterization of $U$ \cite[Lemma~2.1]{darvas-survey} gives
$S\leq U$ and $U_\delta\leq U$. To obtain the reverse bound up to a
controlled error, let $s=\sqrt{\delta/\alpha}\in(0,1]$ and put
\[
  W_s=(1-s)U+sS.
\]
This continuous $\theta$-plurisubharmonic function has the same boundary
values as $U_\delta$. Positivity of mixed Bedford--Taylor products,
obtained by smooth approximation and
\cite[Proposition~2.3]{bedford-taylor}, gives
\[
  (\theta+\ddc W_s)^2
    =\bigl((1-s)\Omega_U+s\eta\bigr)^2
    \geq s^2\eta^2
    \geq\delta\chi^2.
\]
Apply the comparison principle~\cite[Theorem~2.3]{blocki} to $W_s-S$ and
$U_\delta-S$ with K\"ahler background $\eta$. The larger Monge--Amp\`ere
measure of the first function and their common boundary values imply
$W_s\leq U_\delta$. Consequently,
\[
  0\leq U-U_\delta
    \leq s(U-S)
    \leq\frac{\sqrt\delta}{\sqrt\alpha}
       \lVert U-S\rVert_{C^0(X\times\overline A)},
\]
which proves~\eqref{eq:uniform-approximation}. The forms $\eta$ and
$\chi$, and hence $\alpha$, were fixed before $\delta$ was chosen.
\end{proof}

\begin{corollary}
\label{cor:fixed-disk-degeneration}
Fix $t\in(0,1)$ and a closed coordinate disk $K\subset C$ centered at
$p$, such that $u_t$ is smooth on a neighborhood of $K$ and
$(\omega+\ddc_Xu_t)(p)=0$. Define the positive smooth function
$\lambda_{\delta,t}$ by
\[
  \omega+\ddc_XU_\delta(\cdot,e^t)
    =\lambda_{\delta,t}\omega.
\]
Then
\begin{equation}
\label{eq:minimum-degeneration}
  \lim_{\delta\downarrow0}\min_K\lambda_{\delta,t}=0.
\end{equation}
In particular, the fiber restrictions in
\eqref{eq:nondegenerate-dirichlet} have no positive lower bound
independent of $\delta$.
\end{corollary}

\begin{proof}
Let $p$ be the center of $K$. The limiting fiber form is smooth on a
neighborhood of $K$ and can be written as
\[
  \omega+\ddc_Xu_t=\lambda_t\omega,
  \qquad \lambda_t(p)=0.
\]
Uniform convergence on this fixed fiber gives distributional convergence
of its complex Hessians. Indeed, for every
$\varphi\in C_c^\infty(K^\circ)$, integration by parts yields
\[
  \int_K\varphi\,\ddc_X
       \bigl(U_\delta(\cdot,e^t)-u_t\bigr)
    =\int_K\bigl(U_\delta(\cdot,e^t)-u_t\bigr)\ddc_X\varphi
    \longrightarrow0.
\]
If~\eqref{eq:minimum-degeneration} failed, positivity of
$\lambda_{\delta,t}$ would give an $\varepsilon>0$ and a sequence
$\delta_k\downarrow0$ with $\lambda_{\delta_k,t}\geq\varepsilon$ on $K$.
Passing to the limit against nonnegative test functions would imply
\[
  \int_K\varphi\lambda_t\omega
    \geq\varepsilon\int_K\varphi\omega
  \qquad
  (0\leq\varphi\in C_c^\infty(K^\circ)).
\]
By continuity and $\lambda_t(p)=0$, there is a disk
$D\Subset K^\circ$ on which $\lambda_t<\varepsilon/2$. A nonzero
nonnegative test function supported in $D$ contradicts this inequality.
\end{proof}

For either constructed pair, the hypotheses hold for every $t\in J$
and $p\in B$ once $K$ is chosen sufficiently small. The endpoints,
time $t$, and disk $K$ are fixed in this argument. For
every $\varepsilon>0$ and all sufficiently small $\delta$, there is a
point $x_\delta\in K$ with $\lambda_{\delta,t}(x_\delta)<\varepsilon$.
Thus a positive lower bound uniform in $\delta$ and valid on all fibers
is impossible, giving the counterexample requested in~\cite[Problem~2.5]{aim}.
The conclusion is about the minimum on $K$; it does not assert convergence
of $\lambda_{\delta,t}$ at the ramification point itself.

Taking the explicit pair on $\PP^1$ and $t\in[1/4,3/4]\subset I$
gives Corollary~\ref{cor:main-approximation} from
Proposition~\ref{prop:approximation-convergence} and
Corollary~\ref{cor:fixed-disk-degeneration}. The same argument applies
to every curve pair of Theorem~\ref{thm:curves} for $t\in(3/8,5/8)$.

\subsection{Products}
\label{subsec:product-examples}

\begin{corollary}
\label{cor:higher-dimensional-products}
Let $(C,\omega_C)$ be any compact connected Riemann surface and choose
the endpoints of Theorem~\ref{thm:curves}, with weak geodesic $U$ and
ramification set $B$. Let $(Y,\omega_Y)$ be any compact K\"ahler
manifold of dimension $n-1$, where $n\geq2$, and set
\[
  M=C\times Y,\qquad
  \omega_M=\pr_C^*\omega_C+\pr_Y^*\omega_Y,\qquad
  u_j^M=u_j\circ\pr_C.
\]
The weak geodesic joining these endpoints is exactly
$U^M(z,y,\tau)=U(z,\tau)$. It is smooth near
$(\{p\}\times Y)\times\{\tau:3/8<\log|\tau|<5/8\}$ for each
$p\in B$. The total form has rank $n$ in these neighborhoods.
For each $t\in(3/8,5/8)$, the fiber form has rank $n-1$ along each
$\{p\}\times Y$ and is positive at nearby points off this hypersurface.
If the ramification index at $p$ is $m\geq2$, the fiber volume ratio
$(\omega_M+\ddc_M u_t^M)^n/\omega_M^n$ vanishes to order exactly
$2m-2$ along $\{p\}\times Y$.

For $n=1$, omit the $Y$ factor: the ranks are $1$ and $0$, and the
same vanishing-order statement holds. For the explicit
$\PP^1$ pair, the same product conclusions hold throughout $t\in I$,
with $B=\{0,\infty\}$ and $m=2$, in particular for
$t\in[1/4,3/4]$.
\end{corollary}

\begin{proof}
For $n\geq2$, the endpoint forms are positive by the split product
formula. On $M\times A$, let $q(z,y,\tau)=(z,\tau)$ and define
$\widetilde U=U\circ q$. Its total
form is
\[
  \widetilde\Omega
    =q^*\Omega_U+\pr_Y^*\omega_Y\geq0,
  \qquad \Omega_U=\pr_C^*\omega_C+\ddc U.
\]
The Bedford--Taylor product satisfies
\begin{equation}
\label{eq:product-hcma}
  \widetilde\Omega^{n+1}
    =\binom{n+1}{2}\,
      q^*(\Omega_U^2)\wedge\pr_Y^*(\omega_Y^{n-1})=0.
\end{equation}
To justify the identity, work in product coordinate domains and
approximate a continuous local plurisubharmonic potential of $\Omega_U$
by a decreasing sequence of smooth plurisubharmonic functions. Since the
limiting potential is continuous, Dini's theorem gives uniform convergence
on compact subsets. In the smooth binomial expansion, terms with at least
three factors from $C\times A$ vanish by dimension, and terms with at
most one such factor contain at least $n$ factors from $Y$ and also vanish.
The only remaining term is the one displayed in~\eqref{eq:product-hcma}.
Continuity of Bedford--Taylor products under locally uniform convergence
\cite[Proposition~2.3]{bedford-taylor} passes this identity to the limit,
and $\Omega_U^2=0$ gives its vanishing.

The function $\widetilde U$ is continuous and bounded, has the prescribed
boundary values, and is rotation-invariant. Uniqueness for bounded
Bedford--Taylor solutions~\cite[Lemma~2.1]{darvas-survey} therefore
identifies it as the weak geodesic joining $u_0^M,u_1^M$.
On the curve model's agreement region the total form $\Omega_U$ has
rank one, so the split sum $\widetilde\Omega$ has rank $n$.
At a ramification point $p$, its fiber restriction is
\[
  0\oplus\omega_Y
  \qquad\text{on }T_pC\oplus TY,
\]
of rank $n-1$. Choose a coordinate $\zeta$ centered at $p$ and write
$\omega_C=\ii h(\zeta)\,d\zeta\wedge d\bar\zeta$, with $h>0$.
In product coordinates the fiber form is
$\ii g_t(\zeta)\,d\zeta\wedge d\bar\zeta+\omega_Y$.
In the top exterior powers of the fiber and background forms,
the common factor $n\,\omega_Y^{n-1}$ cancels, giving
\[
  \frac{(\omega_M+\ddc_M u_t^M)^n}{\omega_M^n}
  =\pr_C^*\!\left(\frac{\omega_C+\ddc_C u_t}{\omega_C}\right)
  =|\zeta|^{2m-2}\frac{A_t(\zeta)}{h(\zeta)} .
\]
Here $g_t=|\zeta|^{2m-2}A_t$ by the factorization preceding
\eqref{eq:curve-vanishing-order}, and $A_t/h$ is smooth and strictly
positive near $p$. Changing the local defining function for
$\{p\}\times Y$ multiplies this last factor by a smooth positive unit,
so the asserted order is intrinsic. The expansion in
\eqref{eq:curve-vanishing-order} is uniform on compact subintervals of
$(3/8,5/8)$.
For $n=1$ these are precisely the curve statements. For the explicit
$\PP^1$ pair, Proposition~\ref{prop:explicit} supplies the agreement
throughout $I$, and~\eqref{eq:expansion} gives the quadratic normal
zeros.
\end{proof}

\section{Discussion}\label{sec:discussion}

The small corrections in Remark~\ref{rem:small} do not place the endpoints in
the small-data regime of Hu's lower-bound theorem
\cite[Theorems~1.1--1.2]{hu}. Their size is measured relative to the
semipositive model, whereas Hu controls the full boundary data relative to a
fixed K\"ahler background. In Hu's Theorem~1.1 the spatial $C^2$ smallness
is uniform over the entire boundary, and the equation parameter must also
be sufficiently small. His Theorem~1.2 instead assumes a small sum of the two
endpoint $C^2$ norms. For the explicit endpoints of
Proposition~\ref{prop:explicit},
\[
 (u_j)_{z\bar z}(0)=-2+10^{-6},\qquad \omega_{z\bar z}(0)=2.
\]
Changing the reference form to $\omega_0=\omega+\ddc u_0$ does not turn the
endpoint difference $f=u_1-u_0$ into small data: at $0$,
\[
 (\omega_0)_{z\bar z}=10^{-6},\qquad df=0,\qquad
 f_{zz}=-2\tanh(1).
\]
Since $df(0)=0$, the last derivative is also the corresponding component of
the covariant Hessian there. Thus its norm relative to $\omega_0$ is large. A change of background
also changes the geometric smallness threshold in Hu's theorem, so the
original threshold cannot be reused.
The distinction between degeneration of the fiber form and the rank of the
total form is recorded in Remark~\ref{rem:rank}.

The standard Hele--Shaw ray recovers our semipositive model, as explained in
Remark~\ref{rem:hele-shaw}. The arc construction of Ross--Witt Nystr\"om has a
different obstruction to choosing strictly K\"ahler endpoints. Passing from
their disk solution to its rotation-invariant weak ray uses a holomorphic
change of variables and the corresponding Fubini--Study potential
correction. At every positive time, this ray has a nonempty spatial open set
on which its fiber form vanishes
\cite[Lemma~2.4 and Propositions~2.6--2.7]{ross-maximal}. Consequently, no
nontrivial finite truncation has two strictly K\"ahler endpoints.

There is also an obstruction to a particular local repair of such an open
zero set. Let $D$ be a relatively compact smooth domain on which the
original fiber form is zero. If a smooth correction $b\geq0$ vanishes on
$\partial D$ and the corrected fiber form remains semipositive on $D$, then
$\ddc b\geq0$ there. The maximum principle gives $b\leq0$, hence $b=0$ on
$D$. This only rules out repairs with the stated boundary values; it does
not rule out changing those values or the underlying Hele--Shaw data. In
our model, the zeros are isolated and the surrounding fiber form is
positive. That positivity absorbs the negative Hessian of the cutoff, while
Theorem~\ref{thm:preservation} preserves the model in the interior. The
corrected geodesic is obtained from the annular envelope problem; its
identification with a Hele--Shaw flow is not needed.

\section*{Use of generative AI}
ChatGPT was used in developing and revising this manuscript, including literature searches (especially finding and modifying the main ansatz), calculations (construction of $b(z)$), examination of the arguments, and language editing. AI assistance does not constitute independent mathematical verification. The authors take full responsibility for any submitted version.

\end{document}